\documentclass[12pt]{article}%
\usepackage{amsfonts}
\usepackage{amsmath}
\usepackage{amssymb}
\usepackage{graphicx}%
\newtheorem{theorem}{Theorem}

\newtheorem{corollary}[theorem]{Corollary}

\newtheorem{definition}[theorem]{Definition}
\newtheorem{example}[theorem]{Example}

\newtheorem{proposition}[theorem]{Proposition}
\newtheorem{remark}[theorem]{Remark}

\newenvironment{proof}[1][Proof]{\noindent\textbf{#1.} }{\ \rule{0.5em}{0.5em}}
\begin{document}

\title{Recent developments in applied pseudoanalytic function theory}
\author{Vladislav V. Kravchenko\\Department of Mathematics\\CINVESTAV, Queretaro\\Libramiento Norponiente No. 2000\\Fracc. Real de Juriquilla\\Queretaro, Qro.\\C.P. 76230\\MEXICO\\e-mail: vkravchenko@qro.cinvestav.mx}
\maketitle

\begin{abstract}
We present recently obtained results in the theory of pseudoanalytic functions
and its applications to elliptic second-order equations. The operator
$(\operatorname{div}p\operatorname{grad}+q)$ with $p$ and $q$ being real
valued functions is factorized with the aid of Vekua type operators of a
special form and as a consequence the elliptic equation
\begin{equation}
(\operatorname{div}p\operatorname{grad}+q)u=0,\label{1}%
\end{equation}
reduces to a homogeneous Vekua equation describing generalized analytic (or
pseudoanalytic) functions.\ As a tool for solving the Vekua equation we use
the theory of Taylor and Laurent series in formal powers for pseudoanalytic
functions developed by L. Bers. The series possess many important properties
of the usual analytic power series. Their applications until recently were
limited mainly because of the impossibility of their explicit construction in
a general situation. We obtain an algorithm which in a really broad range of
practical applications allows us to construct the formal powers and hence the
pseudoanalytic Taylor series in explicit form precisely for the Vekua equation
related to equation (1). In other words, in a bounded domain this gives us a
complete (in $C$-norm) system of exact solutions of (1).

\end{abstract}

\section{Introduction}

The foundations of pseudoanalytic function theory have been created by a
considerable number of mathematicians among which Lipman Bers and Ilya Vekua
played the most prominent role. In the works of I. Vekua and many other
researchers pseudoanalytic functions are called generalized analytic.
Nevertheless in the present work we use the term \textquotedblleft
pseudoanalytic\textquotedblright\ in order to emphasize the fact that we
mainly use the part of the theory developed by L. Bers and his collaborators.

In the recent works of the author \cite{KrBers}, \cite{Krpseudoan} and
\cite{KrJPhys06} a close connection between the second-order elliptic
equation
\begin{equation}
(\operatorname{div}p\operatorname{grad}+q)u=0 \label{maineqIntro}%
\end{equation}
and a Vekua equation of a special form was presented. This connection is a
direct generalization of a relation which exists between harmonic and analytic
functions. The special form of the arising Vekua equation (which we call the
main Vekua equation) allows us to apply well developed methods of
pseudoanalytic function theory (\cite{Begehr}, \cite{Berskniga},
\cite{Courant}, \cite{Tutschke}, \cite{Vekua}, \cite{Wendland} and others) and
of $p$-analytic function theory \cite{Polozhy} to the analysis of the
corresponding second-order equations. In this work we restrict ourselves to
the development of the theory of series in formal powers for the main Vekua
equation and as a consequence of the corresponding theory for equation
(\ref{maineqIntro}).

Formal powers were defined by L. Bers (see \cite{Berskniga}) and represent a
generalization of the usual powers $(z-z_{0})^{n}$ which play a crucial role
in the one-dimensional complex analysis. As their name reveals formal powers
in general are not powers. They behave as $(z-z_{0})^{n}$ only locally, near
the center, and in fact can be complex functions of a quite arbitrary nature.
Nevertheless they are solutions of a corresponding Vekua equation and under
quite general conditions represent a complete system of its solutions in the
same sense as any analytic function under quite general and well known
conditions can be represented as its normally convergent Taylor series.

The main result of the present work is a procedure for explicit construction
of formal powers corresponding to the main Vekua equation in a very general
situation. From the relation of the main Vekua equation to equation
(\ref{maineqIntro}) we obtain that under quite general conditions we are able
to construct explicitly a complete system of solutions of (\ref{maineqIntro}).
More precisely, let us consider, e.g., the conductivity equation%
\[
\operatorname{div}\left(  p\operatorname{grad}u\right)  =0.
\]
Our result then gives us the possibility to construct explicitly a complete
system of solutions of this equation if $p$ has the form%
\begin{equation}
p=\Phi(\varphi)\Psi(\psi) \label{sepform}%
\end{equation}
where $(\varphi,\psi)$ is any orthogonal coordinate system, $\Phi$ and $\Psi$
are arbitrary positive differentiable functions.

In the case of the stationary Schr\"{o}dinger equation%
\begin{equation}
(-\Delta+q)u=0 \label{SchrIntro}%
\end{equation}
in order to construct a complete system of solutions explicitly we need a
particular solution of this equation of the form (\ref{sepform}). Note that
before this result has been obtained the knowledge of one particular solution
of a second-order equation in two dimensions like (\ref{SchrIntro}) had not
given much information about the general solution. Now one particular solution
of (\ref{SchrIntro}) is a generator of a complete system of solutions of
(\ref{SchrIntro}) which in a sense and for many purposes represents the
general solution of the equation.

The present work is an introduction to this new method and includes
explanation of the theory behind it and some first examples of application.
The paper is organized as follows. In section 2 we present the already known
results from pseudoanalytic function theory concerning mainly the theory of
series in formal powers. In section 3 we give our recent results on the
relationship between equation (\ref{maineqIntro}) and the main Vekua equation.
In section 4 we explain how a complete system of solutions of equation
(\ref{maineqIntro}) is obtained from a complete system of pseudoanalytic
formal powers. Section 5 is dedicated to a classical but unfortunately not
widely known relation between orthogonal coordinate systems in the plane and
complex analytic functions. In section 6 we first explain and then prove the
main new result of this work which opens the way for explicit construction of
formal powers for the main Vekua equation. In section 7 we explain how this
construction gives us complete systems of solutions of second-order elliptic
equations. We present some examples which were calculated using the Matlab
tool for symbolic calculation. The choice of the examples was inspired by the
desire to show another interesting relation between our work and the theory of
so called $\nu$-regular functions (solutions of Vekua equations with constant
coefficients) started by R. J. Duffin in \cite{Duffin71} and \cite{Duffin72}.
We show that the complete systems of solutions of Yukawa and Helmholtz
equations obtained in \cite{BraggDettman95} represent an example of a complete
system of solutions for a second-order elliptic equation constructed on the
base of pseudoanalytic formal powers.

Note that some of the results presented in this work are sufficiently simple
that may serve as a nice and useful addition to a standard complex analysis
course, for example, the construction of conjugate metaharmonic functions (see
subsection 3.3). The author has already done this didactical experiment
successfully. At the same time it is clear that very close to these simple
things there are many open and not so simple problems. First of all, the
theory of pseudoanalytic (or generalized analytic) functions is very rich and
meanwhile some of the well known results from this theory like, e.g., the
Cauchy integral theorem and the Morera theorem have been already transferred
to the corresponding second-order elliptic equation (see \cite{KrBers}), many
others still wait for their application. For example, it is of great interest
and importance to construct not only positive formal powers (as we do in the
present work) but also the negative formal powers for the main Vekua equation
and hence the Cauchy kernel and integral representations for solutions of
second-order equations. Moreover, in the present work we consider equation
(\ref{maineqIntro}) with real valued coefficients. Nevertheless as was shown
in \cite{KrJPhys06} the same scheme works in the case of complex valued
coefficients in (\ref{maineqIntro}), but the corresponding main Vekua equation
is bicomplex. The theory of bicomplex Vekua equations (which are closely
related to the Dirac equation with electromagnetic potential) we started to
develop in \cite{KrAntonio}, but it is clear that it is much more complicated
than the theory of complex Vekua equations and represents a separate and
important challenge.

The development in three or more dimensions of a theory comparable with the
theory of pseudoanalytic functions in the plane was an object of study of many
researchers (see, e.g., \cite{Mal98}). In general, the results obtained in
this direction are much less complete than their two-dimensional counterparts.
In \cite{KrJPhys06} a quaternionic generalization of the main Vekua equation
was obtained, an equation which clearly preserves the basic important
properties which made the present work possible. It is interesting to develop
the L. Bers theory for this main quaternionic Vekua equation which seems to be
possible due to its quite special form.

Finally, behind the results presented here there is an essential and in a
sense universal idea of factorization. The main Vekua equation (as we show in
section 3) is a result of a factorization of a corresponding second-order
operator. This kind of factorization is studied not only in the case of
elliptic operators (\cite{Swansolo}, \cite{Swan}, \cite{KKW}, \cite{AQA},
\cite{KrBers}, \cite{KrJPhys06}) but as well in the case of operators of other
types (\cite{BernsteinCV2006} and \cite{CKS2005}). Thus, it seems quite
possible that similar results can be obtained, e.g., not only for a stationary
Schr\"{o}dinger equation (as in the present work) but also for a
time-dependent Schr\"{o}dinger equation.

\section{Some definitions and results from pseudoanalytic function theory}

This section is based on notions and results presented in \cite{Berskniga} and
\cite{BersStat}. Let $\Omega$ be a domain in $\mathbf{R}^{2}$. Throughout the
whole paper we suppose that $\Omega$ is a simply connected domain.

\subsection{Generating pair, derivative and antiderivative}

\begin{definition}
A pair of complex functions $F$ and $G$ possessing in $\Omega$ partial
derivatives with respect to the real variables $x$ and $y$ is said to be a
generating pair if it satisfies the inequality
\begin{equation}
\operatorname{Im}(\overline{F}G)>0\qquad\text{in }\Omega. \label{GenPairCond}%
\end{equation}

\end{definition}

Denote $\partial_{\overline{z}}=\frac{1}{2}\left(  \frac{\partial}{\partial
x}+i\frac{\partial}{\partial y}\right)  $ and $\partial_{z}=\frac{1}{2}\left(
\frac{\partial}{\partial x}-i\frac{\partial}{\partial y}\right)  $. The
following expressions are known as characteristic coefficients of the pair
$(F,G)$
\[
a_{(F,G)}=-\frac{\overline{F}G_{\overline{z}}-F_{\overline{z}}\overline{G}%
}{F\overline{G}-\overline{F}G},\qquad b_{(F,G)}=\frac{FG_{\overline{z}%
}-F_{\overline{z}}G}{F\overline{G}-\overline{F}G},
\]

\[
A_{(F,G)}=-\frac{\overline{F}G_{z}-F_{z}\overline{G}}{F\overline{G}%
-\overline{F}G},\qquad B_{(F,G)}=\frac{FG_{z}-F_{z}G}{F\overline{G}%
-\overline{F}G},
\]
where the subindex $\overline{z}$ or $z$ means the application of
$\partial_{\overline{z}}$ or $\partial_{z}$ respectively.

Every complex function $W$ defined in a subdomain of $\Omega$ admits the
unique representation $W=\phi F+\psi G$ where the functions $\phi$ and $\psi$
are real valued. Thus, the pair $(F,G)$ generalizes the pair $(1,i)$ which
corresponds to usual complex analytic function theory. 

Sometimes it is convenient to associate with the function $W$ the function
$\omega=\phi+i\psi$. The correspondence between $W$ and $\omega$ is one-to-one.

The $(F,G)$-derivative $\overset{\cdot}{W}=\frac{d_{(F,G)}W}{dz}$ of a
function $W$ exists and has the form
\begin{equation}
\overset{\cdot}{W}=\phi_{z}F+\psi_{z}G=W_{z}-A_{(F,G)}W-B_{(F,G)}\overline{W}
\label{FGder}%
\end{equation}
if and only if
\begin{equation}
\phi_{\overline{z}}F+\psi_{\overline{z}}G=0. \label{phiFpsiG}%
\end{equation}
This last equation can be rewritten also in the following form%
\[
W_{\overline{z}}=a_{(F,G)}W+b_{(F,G)}\overline{W}%
\]
which we call a Vekua equation. Solutions of this equation are called
$(F,G)$-pseudoanalytic functions. If $W$ is $(F,G)$-pseudoanalytic, the
associated function $\omega$ is called $(F,G)$-pseudoanalytic of second kind.

\begin{remark}
\label{RemDerivativeOfGenerators}The functions $F$ and $G$ are $(F,G)$%
-pseudoanalytic, and $\overset{\cdot}{F}\equiv\overset{\cdot}{G}\equiv0$.
\end{remark}

\begin{definition}
\label{DefSuccessor}Let $(F,G)$ and $(F_{1},G_{1})$ - be two generating pairs
in $\Omega$. $(F_{1},G_{1})$ is called \ successor of $(F,G)$ and $(F,G)$ is
called predecessor of $(F_{1},G_{1})$ if%
\[
a_{(F_{1},G_{1})}=a_{(F,G)}\qquad\text{and}\qquad b_{(F_{1},G_{1})}%
=-B_{(F,G)}\text{.}%
\]

\end{definition}

The importance of this definition becomes obvious from the following statement.

\begin{theorem}
\label{ThBersDer}Let $W$ be an $(F,G)$-pseudoanalytic function and let
$(F_{1},G_{1})$ be a successor of $(F,G)$. Then $\overset{\cdot}{W}$ is an
$(F_{1},G_{1})$-pseudoanalytic function.
\end{theorem}

This theorem shows us that to the difference of analytic functions whose
derivatives are again analytic, the $(F,G)$-derivatives of pseudoanalytic
functions are in general solutions of another Vekua equation.

\begin{definition}
\label{DefAdjoint}Let $(F,G)$ be a generating pair. Its adjoint generating
pair $(F,G)^{\ast}=(F^{\ast},G^{\ast})$ is defined by the formulas%
\[
F^{\ast}=-\frac{2\overline{F}}{F\overline{G}-\overline{F}G},\qquad G^{\ast
}=\frac{2\overline{G}}{F\overline{G}-\overline{F}G}.
\]

\end{definition}

The $(F,G)$-integral is defined as follows
\[
\int_{\Gamma}Wd_{(F,G)}z=F(z_{1})\operatorname{Re}\int_{\Gamma}G^{\ast
}Wdz+G(z_{1})\operatorname{Re}\int_{\Gamma}F^{\ast}Wdz
\]
where $\Gamma$ is a rectifiable curve leading from $z_{0}$ to $z_{1}$.

If $W=\phi F+\psi G$ is an $(F,G)$-pseudoanalytic function where $\phi$ and
$\psi$ are real valued functions then
\begin{equation}
\int_{z_{0}}^{z}\overset{\cdot}{W}d_{(F,G)}z=W(z)-\phi(z_{0})F(z)-\psi
(z_{0})G(z), \label{FGAnt}%
\end{equation}
and as $\overset{\cdot}{F}=\overset{}{\overset{\cdot}{G}=}0$, this integral is
path-independent and represents the $(F,G)$-antiderivative of $\overset{\cdot
}{W}$.

A continuous function $w$ defined in a domain $\Omega$ is called
$(F,G)$-integrable if for every closed curve $\Gamma$ situated in a simply
connected subdomain of $\Omega$,
\[
\int_{\Gamma}wd_{(F,G)}z=0
\]
or what is the same%
\begin{equation}
\operatorname{Re}\int_{\Gamma}G^{\ast}wdz+i\operatorname{Re}\int_{\Gamma
}F^{\ast}wdz=0. \label{Cauchyth}%
\end{equation}

\begin{theorem}
\label{ThCauchyPa}An $(F,G)$-derivative $\overset{\cdot}{W}$ of an
$(F,G)$-pseudoanalytic function $W$ is $(F,G)$-integrable.
\end{theorem}

\begin{theorem}
\label{ThMoreraPa}Let $(F,G)$ be a predecessor of $(F_{1},G_{1})$. A
continuous function is $(F_{1},G_{1})$-pseudoanalytic if and only if it is
$(F,G)$-integrable.

\begin{remark}
It is easy to see that in the case $F\equiv1$, $G\equiv i$, equality
(\ref{Cauchyth}) turns into the Cauchy integral theorem for analytic
functions: $\int_{\Gamma}wdz=0$.
\end{remark}
\end{theorem}

\subsection{Generating sequences and Taylor series in formal
powers\label{SubsectGenSeq}}

In order to introduce the notion of pseudoanalytic derivatives of arbitrary
order the following definition is necessary.

\begin{definition}
\label{DefSeq}A sequence of generating pairs $\left\{  (F_{m},G_{m})\right\}
$, $m=0,\pm1,\pm2,\ldots$ , is called a generating sequence if $(F_{m+1}%
,G_{m+1})$ is a successor of $(F_{m},G_{m})$. If $(F_{0},G_{0})=(F,G)$, we say
that $(F,G)$ is embedded in $\left\{  (F_{m},G_{m})\right\}  $.
\end{definition}

\begin{theorem}
Let \ $(F,G)$ be a generating pair in $\Omega$. Let $\Omega_{1}$ be a bounded
domain, $\overline{\Omega}_{1}\subset\Omega$. Then $(F,G)$ can be embedded in
a generating sequence in $\Omega_{1}$.
\end{theorem}

\begin{definition}
A generating sequence $\left\{  (F_{m},G_{m})\right\}  $ is said to have
period $\mu>0$ if $(F_{m+\mu},G_{m+\mu})$ is equivalent to $(F_{m},G_{m})$
that is their characteristic coefficients coincide.
\end{definition}

Let $W$ be an $(F,G)$-pseudoanalytic function. Using a generating sequence in
which $(F,G)$ is embedded we can define the higher derivatives of $W$ by the
recursion formula%
\[
W^{[0]}=W;\qquad W^{[m+1]}=\frac{d_{(F_{m},G_{m})}W^{[m]}}{dz},\quad
m=1,2,\ldots\text{.}%
\]

A generating sequence defines an infinite sequence of Vekua equations. If for
a given (original) Vekua equation we know not only a corresponding generating
pair but the whole generating sequence, that is a couple of exact and
independent solutions for each of the Vekua equations from the infinite
sequence of equations corresponding to the original one, we are able to
construct an infinite system of solutions of the original Vekua equation as is
shown in the next definition. Moreover, as we show in the next subsection
under quite general conditions this infinite system of solutions is complete.

\begin{definition}
\label{DefFormalPower}The formal power $Z_{m}^{(0)}(a,z_{0};z)$ with center at
$z_{0}\in\Omega$, coefficient $a$ and exponent $0$ is defined as the linear
combination of the generators $F_{m}$, $G_{m}$ with real constant coefficients
$\lambda$, $\mu$ chosen so that $\lambda F_{m}(z_{0})+\mu G_{m}(z_{0})=a$. The
formal powers with exponents $n=1,2,\ldots$ are defined by the recursion
formula%
\begin{equation}
Z_{m}^{(n+1)}(a,z_{0};z)=(n+1)\int_{z_{0}}^{z}Z_{m+1}^{(n)}(a,z_{0}%
;\zeta)d_{(F_{m},G_{m})}\zeta. \label{recformula}%
\end{equation}

\end{definition}

This definition implies the following properties.

\begin{enumerate}
\item $Z_{m}^{(n)}(a,z_{0};z)$ is an $(F_{m},G_{m})$-pseudoanalytic function
of $z$.

\item If $a^{\prime}$ and $a^{\prime\prime}$ are real constants, then
$Z_{m}^{(n)}(a^{\prime}+ia^{\prime\prime},z_{0};z)=a^{\prime}Z_{m}%
^{(n)}(1,z_{0};z)+a^{\prime\prime}Z_{m}^{(n)}(i,z_{0};z).$

\item The formal powers satisfy the differential relations%
\[
\frac{d_{(F_{m},G_{m})}Z_{m}^{(n)}(a,z_{0};z)}{dz}=nZ_{m+1}^{(n-1)}%
(a,z_{0};z).
\]

\item The asymptotic formulas
\begin{equation}
Z_{m}^{(n)}(a,z_{0};z)\sim a(z-z_{0})^{n},\quad z\rightarrow z_{0}
\label{asymptformulas}%
\end{equation}
hold.
\end{enumerate}

Assume now that
\begin{equation}
W(z)=\sum_{n=0}^{\infty}Z^{(n)}(a,z_{0};z) \label{series}%
\end{equation}
where the absence of the subindex $m$ means that all the formal powers
correspond to the same generating pair $(F,G),$ and the series converges
uniformly in some neighborhood of $z_{0}$. It can be shown that the uniform
limit of pseudoanalytic functions is pseudoanalytic, and that a uniformly
convergent series of $(F,G)$-pseudoanalytic functions can be $(F,G)$%
-differentiated term by term. Hence the function $W$ in (\ref{series}) is
$(F,G)$-pseudoanalytic and its $r$th derivative admits the expansion
\[
W^{[r]}(z)=\sum_{n=r}^{\infty}n(n-1)\cdots(n-r+1)Z_{r}^{(n-r)}(a_{n}%
,z_{0};z).
\]
From this the Taylor formulas for the coefficients are obtained%
\begin{equation}
a_{n}=\frac{W^{[n]}(z_{0})}{n!}. \label{Taylorcoef}%
\end{equation}

\begin{definition}
Let $W(z)$ be a given $(F,G)$-pseudoanalytic function defined for small values
of $\left\vert z-z_{0}\right\vert $. The series%
\begin{equation}
\sum_{n=0}^{\infty}Z^{(n)}(a,z_{0};z) \label{Taylorseries}%
\end{equation}
with the coefficients given by (\ref{Taylorcoef}) is called the Taylor series
of $W$ at $z_{0}$, formed with formal powers.
\end{definition}

The Taylor series always represents the function asymptotically:%
\begin{equation}
W(z)-\sum_{n=0}^{N}Z^{(n)}(a,z_{0};z)=O\left(  \left\vert z-z_{0}\right\vert
^{N+1}\right)  ,\quad z\rightarrow z_{0}, \label{asympt}%
\end{equation}
for all $N$. This implies (since a pseudoanalytic function can not have a zero
of arbitrarily high order without vanishing identically) that the sequence of
derivatives $\left\{  W^{[n]}(z_{0})\right\}  $ determines the function $W$ uniquely.

If the series (\ref{Taylorseries}) converges uniformly in a neighborhood of
$z_{0}$, it converges to the function $W$.

\subsection{Convergence theorems\label{SubsectConvergence}}

The statements given in this subsection were obtained by L. Bers
\cite{Berskniga}, \cite{BersFormalPowers} and S. Agmon and L. Bers
\cite{AgmonBers}.

\begin{theorem}
\label{ThConvPer}\cite{Berskniga} The formal Taylor expansion
(\ref{Taylorseries}) of a pseudoanalytic function in formal powers defined by
a periodic generating sequence converges in some neighborhood of the center.
\end{theorem}

This theorem means only a local completeness of the system of formal powers.
The following definition due to L. Bers describes the case when corresponding
formal powers represent a globally complete system of solutions of a Vekua
equation much as in the case of usual powers of the variable $z$ and the
Cauchy-Riemann equation.

\begin{definition}
\cite{Berskniga} A generating pair $(F,G)$ is called complete if these
functions are defined and satisfy the H\"{o}lder condition for all finite
values of $z$, the limits $F(\infty)$, $G(\infty)$ exist, $\operatorname{Im}%
(\overline{F(\infty)}G(\infty))>0$, and the functions $F(1/z)$, $G(1/z)$ also
satisfy the H\"{o}lder condition. \ A complete generating pair is called
normalized if $F(\infty)=1$, $G(\infty)=i$.
\end{definition}

A generating pair equivalent to a complete one is complete, and every complete
generating pair is equivalent to a uniquely determined normalized pair. The
adjoint of a complete (normalized) generating pair is complete (normalized).

From now on we assume that $(F,G)$ is a complete normalized generating pair.
Then much more can be said on the series of corresponding formal powers. We
limit ourselves to the following completeness results (the expansion theorem
and Runge%
\'{}%
s approximation theorem for pseudoanalytic functions).

Following \cite{Berskniga} we shall say that a sequence of functions $W_{n}$
converges normally in a domain $\Omega$ if it converges uniformly on every
bounded closed subdomain of $\Omega$.

\begin{theorem}
\label{ThTaylorRepr}Let $W$ be an $(F,G)$-pseudoanalytic function defined for
$\left\vert z-z_{0}\right\vert <R$. Then it admits a unique expansion of the
form $W(z)=\sum_{n=0}^{\infty}Z^{(n)}(a_{n},z_{0};z)$ which converges normally
for $\left\vert z-z_{0}\right\vert <\theta R$, where $\theta$ is a positive
constant depending on the generating sequence.
\end{theorem}

The first version of this theorem was proved in \cite{AgmonBers}. We follow
here \cite{BersFormalPowers}.

\begin{remark}
\label{RemNecConditions}Necessary and sufficient conditions for the relation
$\theta=1$ are, unfortunately, not known. However, in \cite{BersFormalPowers}
the following sufficient conditions for the case when the generators $(F,G)$
possess partial derivatives are given. One such condition reads:%
\[
\left\vert F_{\overline{z}}(z)\right\vert +\left\vert G_{\overline{z}%
}(z)\right\vert \leq\frac{\operatorname*{Const}}{1+\left\vert z\right\vert
^{1+\varepsilon}}%
\]
for some $\varepsilon>0$. Another condition is
\[
\int\int_{\left\vert z\right\vert <\infty}\left(  \left\vert F_{\overline{z}%
}\right\vert ^{2-\varepsilon}+\left\vert F_{\overline{z}}\right\vert
^{2+\varepsilon}+\left\vert G_{\overline{z}}\right\vert ^{2-\varepsilon
}+\left\vert G_{\overline{z}}\right\vert ^{2+\varepsilon}\right)  dxdy<\infty
\]
for some $0<\varepsilon<1$.
\end{remark}

\begin{theorem}
\cite{BersFormalPowers}\label{ThRunge} A pseudoanalytic function defined in a
simply connected domain can be expanded into a normally convergent series of
formal polynomials (linear combinations of formal powers with positive exponents).
\end{theorem}

\begin{remark}
This theorem admits a direct generalization onto the case of a multiply
connected domain (see \cite{BersFormalPowers}).
\end{remark}

In posterior works \cite{IsmTagieva}, \cite{Menke}, \cite{Fryant} and others
deep results on interpolation and on the degree of approximation by
pseudopolynomials were obtained. For example,

\begin{theorem}
\cite{Menke}\label{ThMenke} Let $W$ be a pseudoanalytic function in a domain
$\Omega$ bounded by a Jordan curve and satisfy the H\"{o}lder condition on
$\partial\Omega$ with the exponent $\alpha$ ($0<\alpha\leq1$). Then for any
$\varepsilon>0$ and any natural $n$ there exists a pseudopolynomial of order
$n$ satisfying the inequality
\[
\left\vert W(z)-P_{n}(z)\right\vert \leq\frac{\operatorname*{Const}}%
{n^{\alpha-\varepsilon}}\qquad\text{for any }z\in\overline{\Omega}%
\]
where the constant does not depend on $n$, but only on $\varepsilon$.
\end{theorem}

\section{Solutions of second order elliptic equations as real components of
complex pseudoanalytic functions\label{SectSchr}}

\subsection{Factorization of the stationary Schr\"{o}dinger operator}

It is well known that if $f_{0}$ is a nonvanishing particular solution of the
one-dimensional stationary Schr\"{o}dinger equation%
\[
\left(  -\frac{d^{2}}{dx^{2}}+\nu(x)\right)  f(x)=0
\]
then the Schr\"{o}dinger operator can be factorized as follows%
\[
\frac{d^{2}}{dx^{2}}-\nu(x)=\left(  \frac{d}{dx}+\frac{f_{0}^{\prime}}{f_{0}%
}\right)  \left(  \frac{d}{dx}-\frac{f_{0}^{\prime}}{f_{0}}\right)  .
\]
We start with a generalization of this result onto a two-dimensional
situation. Consider the equation%
\begin{equation}
\left(  -\Delta+\nu\right)  f=0 \label{Schrod}%
\end{equation}
in some domain $\Omega\subset\mathbf{R}^{2}$, where $\Delta=\frac{\partial
^{2}}{\partial x^{2}}+\frac{\partial^{2}}{\partial y^{2}}$, $\nu$ and $f$ are
real valued functions. We assume that $f$ is a twice continuously
differentiable function. By $C$ we denote the complex conjugation operator.

\begin{theorem}
\cite{Krpseudoan} \label{Thfact}Let $f$ be a positive in $\Omega$ particular
solution of (\ref{Schrod}). Then for any real valued function $\varphi\in
C^{2}(\Omega)$ the following equalities hold%
\begin{equation}
\frac{1}{4}\left(  \Delta-\nu\right)  \varphi=\left(  \partial_{\overline{z}%
}+\frac{f_{z}}{f}C\right)  \left(  \partial_{z}-\frac{f_{z}}{f}C\right)
\varphi=\left(  \partial_{z}+\frac{f_{\overline{z}}}{f}C\right)  \left(
\partial_{\overline{z}}-\frac{f_{\overline{z}}}{f}C\right)  \varphi.
\label{fact}%
\end{equation}

\end{theorem}

\begin{proof}
Consider%
\begin{align}
\left(  \partial_{\overline{z}}+\frac{f_{z}}{f}C\right)  \left(  \partial
_{z}-\frac{f_{z}}{f}C\right)  \varphi &  =\frac{1}{4}\Delta\varphi
-\frac{\left\vert \partial_{z}f\right\vert ^{2}}{f^{2}}\varphi-\partial
_{\overline{z}}\left(  \frac{\partial_{z}f}{f}\right)  \varphi\nonumber\\
&  =\frac{1}{4}(\Delta\varphi-\frac{\Delta f}{f}\varphi)=\frac{1}{4}\left(
\Delta-\nu\right)  \varphi. \label{vsp2}%
\end{align}
Thus, we have the first equality in (\ref{fact}). Now application of $C$ to
both sides of (\ref{vsp2}) gives us the second equality in (\ref{fact}).
\end{proof}

The operator $\partial_{z}-\frac{f_{z}}{f}I$, where $I$ is the identity
operator, can be represented in the form%
\[
\partial_{z}-\frac{f_{z}}{f}I=f\partial_{z}f^{-1}I.
\]
Let us introduce the notation $P=f\partial_{z}f^{-1}I$. Due to Theorem
\ref{Thfact}, if $f$ is a positive solution of (\ref{Schrod}), the operator
$P$ transforms real valued solutions of (\ref{Schrod}) into solutions of the
Vekua equation
\begin{equation}
\left(  \partial_{\overline{z}}+\frac{f_{z}}{f}C\right)  w=0. \label{Vekua}%
\end{equation}

Note that the operator $\partial_{z}$ applied to a real valued function
$\varphi$ can be regarded as a kind of gradient, and if we know that
$\partial_{z}\varphi=\Phi$ in a whole complex plane or in a convex domain,
where $\Phi=\Phi_{1}+i\Phi_{2}$ is a given complex valued function such that
its real part $\Phi_{1}$ and imaginary part $\Phi_{2}$ satisfy the equation
\begin{equation}
\partial_{y}\Phi_{1}+\partial_{x}\Phi_{2}=0, \label{casirot}%
\end{equation}
then we can reconstruct $\varphi$ up to an arbitrary real constant $c$ in the
following way%
\begin{equation}
\varphi(x,y)=2\left(  \int_{x_{0}}^{x}\Phi_{1}(\eta,y)d\eta-\int_{y_{0}}%
^{y}\Phi_{2}(x_{0},\xi)d\xi\right)  +c \label{Antigr}%
\end{equation}
where $(x_{0},y_{0})$ is an arbitrary fixed point in the domain of interest.

By $A$ we denote the integral operator in (\ref{Antigr}):%
\[
A[\Phi](x,y)=2\left(  \int_{x_{0}}^{x}\Phi_{1}(\eta,y)d\eta-\int_{y_{0}}%
^{y}\Phi_{2}(x_{0},\xi)d\xi\right)  +c.
\]
Note that formula (\ref{Antigr}) can be easily extended to any simply
connected domain by considering the integral along an arbitrary rectifiable
curve $\Gamma$ leading from $(x_{0},y_{0})$ to $(x,y)$%
\[
\varphi(x,y)=2\left(  \int_{\Gamma}\Phi_{1}dx-\Phi_{2}dy\right)  +c.
\]
Thus if $\Phi$ satisfies (\ref{casirot}), there exists a family of real valued
functions $\varphi$ such that $\partial_{z}\varphi=\Phi$, given by the formula
$\varphi=A[\Phi]$.

In a similar way we define the operator $\overline{A}$ corresponding to
$\partial_{\overline{z}}$:%
\[
\overline{A}[\Phi](x,y)=2\left(  \int_{x_{0}}^{x}\Phi_{1}(\eta,y)d\eta
+\int_{y_{0}}^{y}\Phi_{2}(x_{0},\xi)d\xi\right)  +c.
\]

Consider the operator $S=fAf^{-1}I$ applicable to any complex valued function
$w$ such that $\Phi=f^{-1}w$ satisfies condition (\ref{casirot}). Then it is
clear that for such $w$ we have that $PSw=w$.

\begin{proposition}
\cite{KrPanalyt} \label{PrPS}Let $f$ be a positive particular solution of
(\ref{Schrod}) and $w$ be a solution of (\ref{Vekua}). Then the real valued
function $g=Sw$ is a solution of (\ref{Schrod}).
\end{proposition}

\begin{proof}
First of all let us check that the function $\Phi=w/f$ satisfies
(\ref{casirot}). Let $u=\operatorname*{Re}w$ and $v=\operatorname*{Im}w$.
Consider%
\begin{equation}
\partial_{y}\Phi_{1}+\partial_{x}\Phi_{2}=\frac{1}{f}\left(  (\partial
_{y}u+\partial_{x}v)-(\frac{\partial_{y}f}{f}u+\frac{\partial_{x}f}%
{f}v)\right)  . \label{vsp4}%
\end{equation}
Note that equation (\ref{Vekua})\ is equivalent to the system
\[
\partial_{x}u-\partial_{y}v=-\frac{\partial_{x}f}{f}u+\frac{\partial_{y}f}%
{f}v,\qquad\partial_{y}u+\partial_{x}v=\frac{\partial_{y}f}{f}u+\frac
{\partial_{x}f}{f}v
\]
from which we obtain that expression (\ref{vsp4}) is zero. Thus the function
$\Phi$ satisfies (\ref{casirot}) and hence the real valued function
$\varphi=A[w/f]\ $is well defined and satisfies the equation $\partial
_{z}\varphi=w/f$.

Consider the expression
\begin{align}
\partial_{\overline{z}}\partial_{z}(Sw)  &  =\partial_{\overline{z}}\left(
\left(  \partial_{z}f\right)  A[\frac{w}{f}]+w\right) \nonumber\\
&  =\left(  \frac{1}{4}\Delta f\right)  A[\frac{w}{f}]+\left(  \partial
_{z}f\right)  \partial_{\overline{z}}A[\frac{w}{f}]-\frac{\partial_{z}f}%
{f}\overline{w}. \label{vsp6}%
\end{align}
For the expression $\partial_{\overline{z}}A[\frac{w}{f}]$ we have
\begin{align}
\partial_{\overline{z}}A[\frac{w}{f}]  &  =\partial_{z}A[\frac{w}%
{f}]+i\partial_{y}A[\frac{w}{f}]\nonumber\\
&  =\frac{w}{f}-2i\frac{v}{f}=\frac{\overline{w}}{f} \label{vsp3}%
\end{align}
where the following observation was used%
\begin{align*}
\partial_{y}A[\frac{u+iv}{f}](x,y)  &  =2\left(  \int_{x_{0}}^{x}\partial
_{y}\left(  \frac{u(\eta,y)}{f(\eta,y)}\right)  d\eta-\frac{v(x_{0}%
,y)}{f(x_{0},y)}\right)  =\\
&  =-2\left(  \int_{x_{0}}^{x}\partial_{\eta}\left(  \frac{v(\eta,y)}%
{f(\eta,y)}\right)  d\eta-\frac{v(x_{0},y)}{f(x_{0},y)}\right)  =-\frac
{2v(x,y)}{f(x,y)}.
\end{align*}
Thus substitution of (\ref{vsp3}) into (\ref{vsp6}) gives us the equality%
\[
\Delta(Sw)=\nu fA[\frac{w}{f}]=\nu Sw.
\]

\end{proof}

\begin{proposition}
\cite{KrPanalyt} \label{PrSP}Let $g$ be a real valued solution of
(\ref{Schrod}). Then
\[
SPg=g+cf
\]
where $c$ is an arbitrary real constant.
\end{proposition}

\begin{proof}
Consider
\[
SPg=fA\partial_{z}[\frac{g}{f}]=f(\frac{g}{f}+c)=g+cf.
\]

\end{proof}

Theorem \ref{Thfact} together with Proposition \ref{PrPS} show us that
equation (\ref{Schrod}) is equivalent to the Vekua equation (\ref{Vekua}) in
the following sense. Every solution of one of these equations can be
transformed into a solution of the other equation and vice versa.

\subsection{Factorization of the operator $\operatorname{div}%
p\operatorname{grad}+q$.}

The following statement is known in a form of a substitution (see, e.g.,
\cite{Nachman}). Here we formulate it as an operational relation.

\begin{proposition}
\label{PropSchr}Let $p$ and $q$ be complex valued functions, $p\in
C^{2}(\Omega)$ and $p\neq0$ in $\Omega$. Then
\begin{equation}
\operatorname{div}p\operatorname{grad}+q=p^{1/2}(\Delta-r)p^{1/2}\text{\qquad
in }\Omega, \label{firstfact}%
\end{equation}
where
\[
r=\frac{\Delta p^{1/2}}{p^{1/2}}-\frac{q}{p}.
\]

\end{proposition}

\begin{proof}
The easily verified relation
\begin{equation}
\operatorname{div}p\operatorname{grad}=p^{1/2}(\Delta-\frac{\Delta p^{1/2}%
}{p^{1/2}})p^{1/2} \label{factUhlm}%
\end{equation}
is well known (see, e.g., \cite{UhlmannDevelopments}). Adding to both sides of
(\ref{factUhlm}) the term $q$ (and representing it on the right-hand side as
$p^{1/2}\left(  q/p\right)  p^{1/2}$) gives us (\ref{firstfact}).
\end{proof}

\begin{theorem}
\cite{KrJPhys06} \label{ThFactGenSchr}Let\ $p$ and $q$ be real valued
functions, $p\in C^{2}(\Omega)$ and $p\neq0$ in $\Omega$, $u_{0}$ be a
positive particular solution of the equation
\begin{equation}
(\operatorname{div}p\operatorname{grad}+q)u=0\text{\qquad in }\Omega\text{.}
\label{maineq}%
\end{equation}
Then for any real valued continuously twice differentiable function $\varphi$
the following equality holds%
\begin{equation}
\frac{1}{4}(\operatorname{div}p\operatorname{grad}+q)\varphi=p^{1/2}\left(
\partial_{z}+\frac{f_{\overline{z}}}{f}C\right)  \left(  \partial
_{\overline{z}}-\frac{f_{\overline{z}}}{f}C\right)  p^{1/2}\varphi,
\label{mainfact}%
\end{equation}
where
\begin{equation}
f=p^{1/2}u_{0}. \label{fandu}%
\end{equation}

\end{theorem}

\begin{proof}
This is based on (\ref{fact}). From (\ref{firstfact}) we have that if $u_{0}$
is a solution of (\ref{maineq}) then the function (\ref{fandu}) is a solution
of the equation%
\begin{equation}
(\Delta-r)f=0. \label{Schrodr}%
\end{equation}
Then combining (\ref{firstfact}) and (\ref{fact}) we obtain (\ref{mainfact}).
\end{proof}

\begin{remark}
According to (\ref{factUhlm}), $\Delta-r=f^{-1}\operatorname{div}%
f^{2}\operatorname{grad}f^{-1}$ where $f$ is a solution of (\ref{Schrodr}).
Then from (\ref{firstfact}) we have%
\begin{equation}
\operatorname{div}p\operatorname{grad}+q=p^{1/2}f^{-1}\operatorname{div}%
f^{2}\operatorname{grad}f^{-1}p^{1/2}. \label{reldivgrad}%
\end{equation}
Taking into account (\ref{fandu}) we obtain
\[
\operatorname{div}p\operatorname{grad}+q=u_{0}^{-1}\operatorname{div}%
pu_{0}^{2}\operatorname{grad}u_{0}^{-1}\text{\qquad in }\Omega\text{.}%
\]

\end{remark}

\begin{remark}
Let $q\equiv0$. Then $u_{0}$ can be chosen as $u_{0}\equiv1$. Hence
(\ref{mainfact}) gives us the equality
\[
\frac{1}{4}\operatorname{div}(p\operatorname{grad}\varphi)=p^{1/2}\left(
\partial_{z}+\frac{\partial_{\overline{z}}p^{1/2}}{p^{1/2}}C\right)  \left(
\partial_{\overline{z}}-\frac{\partial_{\overline{z}}p^{1/2}}{p^{1/2}%
}C\right)  (p^{1/2}\varphi).
\]

\end{remark}

In what follows we suppose that in $\Omega$ there exists a positive particular
solution of (\ref{maineq}) which we denote by $u_{0}$.

Let $f$ be a real function of $x$ and $y$. Consider the Vekua equation
\begin{equation}
W_{\overline{z}}=\frac{f_{\overline{z}}}{f}\overline{W}\text{\qquad in }%
\Omega\text{.}\label{Vekuamain}%
\end{equation}
This equation plays a crucial role in all that follows hence we will call it
the \textbf{main Vekua equation}.

Denote $W_{1}=\operatorname*{Re}W$ and $W_{2}=\operatorname{Im}W.$

\begin{remark}
\cite{Krpseudoan} \label{RemAnotherForm} Equation (\ref{Vekuamain}) can be
written as follows%
\begin{equation}
f\partial_{\overline{z}}(f^{-1}W_{1})+if^{-1}\partial_{\overline{z}}%
(fW_{2})=0. \label{anotherformVekua}%
\end{equation}

\end{remark}

\begin{theorem}
\cite{KrJPhys06} \label{ThConjugate} Let $W=W_{1}+iW_{2}$ be a solution of
(\ref{Vekuamain}). Then $U=f^{-1}W_{1}$ is a solution of the conductivity
equation
\begin{equation}
\operatorname*{div}(f^{2}\nabla U)=0\qquad\text{in }\Omega, \label{divf2}%
\end{equation}
and $V=fW_{2}$ is a solution of the associated conductivity equation
\begin{equation}
\operatorname*{div}(f^{-2}\nabla V)=0\qquad\text{in }\Omega, \label{divf-2}%
\end{equation}
the function $W_{1}$ is a solution of the stationary Schr\"{o}dinger equation
\begin{equation}
-\Delta W_{1}+r_{1}W_{1}=0\qquad\text{in }\Omega\label{Schr1}%
\end{equation}
with $r_{1}=\Delta f/f,$ and $W_{2}$ is a solution of the associated
stationary Schr\"{o}dinger equation
\begin{equation}
-\Delta W_{2}+r_{2}W_{2}=0\qquad\text{in }\Omega\label{Schr2}%
\end{equation}
where $r_{2}=2(\nabla f)^{2}/f^{2}-r_{1}$ and $(\nabla f)^{2}=f_{x}^{2}%
+f_{y}^{2}$.
\end{theorem}

\begin{proof}
To prove the first part of the theorem we use the form of equation
(\ref{Vekuamain}) given in Remark \ref{RemAnotherForm}. Multiplying
(\ref{anotherformVekua}) by $f$ and applying $\partial_{z}$ gives%
\[
\partial_{z}\left(  f^{2}\partial_{\overline{z}}\left(  f^{-1}W_{1}\right)
\right)  +\frac{i}{4}\Delta\left(  fW_{2}\right)  =0
\]
from where we have that $\operatorname*{Re}\left(  \partial_{z}\left(
f^{2}\partial_{\overline{z}}\left(  f^{-1}W_{1}\right)  \right)  \right)  =0$
which is equivalent to (\ref{divf2}) where $U=f^{-1}W_{1}$.

Multiplying (\ref{anotherformVekua}) by $f^{-1}$ and applying $\partial_{z}$
gives%
\[
\frac{1}{4}\Delta\left(  f^{-1}W_{1}\right)  +i\partial_{z}\left(
f^{-2}\partial_{\overline{z}}\left(  fW_{2}\right)  \right)  =0
\]
from where we have that $\operatorname*{Re}\left(  \partial_{z}\left(
f^{-2}\partial_{\overline{z}}\left(  fW_{2}\right)  \right)  \right)  =0$
which is equivalent to (\ref{divf-2}) where $V=fW_{2}$.

From (\ref{factUhlm}) we have
\[
\left(  \Delta-r_{1}\right)  W_{1}=f^{-1}\operatorname*{div}(f^{2}%
\nabla\left(  f^{-1}W_{1}\right)  ).
\]
Hence from the just proven equation (\ref{divf2}) we obtain that $W_{1}$ is a
solution of (\ref{Schr1}).

In order to obtain equation (\ref{Schr2}) for $W_{2}$ it should be noticed
that
\[
f\operatorname*{div}(f^{-2}\nabla(fW_{2}))=\left(  \Delta-r_{2}\right)  W_{2}.
\]

\end{proof}

\begin{remark}
\label{RelAstalaPaivarinta} Observe that the pair of functions
\begin{equation}
F=f\quad\text{and\quad}G=\frac{i}{f} \label{genpair}%
\end{equation}
is a generating pair for (\ref{Vekuamain}). This allows us to rewrite
(\ref{Vekuamain}) in the form of an equation for pseudoanalytic functions of
second kind (equation (\ref{phiFpsiG}))%
\begin{equation}
\varphi_{\overline{z}}f+\psi_{\overline{z}}\frac{i}{f}=0,
\label{Vekuamain2kind}%
\end{equation}
where $\varphi$ and $\psi$ are real valued functions. If $\varphi$ and $\psi$
satisfy (\ref{Vekuamain2kind}) then $W=\varphi f+\psi\frac{i}{f}$ is a
solution of (\ref{Vekuamain}) and vice versa.

Denote $w=\varphi+\psi i$. Then from (\ref{Vekuamain2kind}) we have%
\[
(w+\overline{w})_{\overline{z}}f+(w-\overline{w})_{\overline{z}}\frac{1}%
{f}=0,
\]
which is equivalent to the equation
\begin{equation}
w_{\overline{z}}=\frac{1-f^{2}}{1+f^{2}}\overline{w}_{\overline{z}}
\label{APform}%
\end{equation}
The relation between (\ref{APform}) and (\ref{divf2}), (\ref{divf-2}) was
observed in \cite{AstalaPaivarinta} and resulted to be essential for solving
the Calder\'{o}n problem in the plane.
\end{remark}

\begin{theorem}
\cite{KrJPhys06}\label{ThConjugate2} Let $W=W_{1}+iW_{2}$ be a solution of
(\ref{Vekuamain}). Assume that $f=p^{1/2}u_{0}$, where $u_{0}$ is a positive
solution of (\ref{maineq}) in $\Omega$. Then $u=p^{-1/2}W_{1}$ is a solution
of (\ref{maineq}) in $\Omega$, and $v=p^{1/2}W_{2}$ is a solution of the
equation
\begin{equation}
(\operatorname*{div}\frac{1}{p}\operatorname*{grad}+q_{1})v=0\qquad\text{in
}\Omega, \label{assocmaineq}%
\end{equation}
where
\begin{equation}
q_{1}=-\frac{1}{p}\left(  \frac{q}{p}+2\left\langle \frac{\nabla p}{p}%
,\frac{\nabla u_{0}}{u_{0}}\right\rangle +2\left(  \frac{\nabla u_{0}}{u_{0}%
}\right)  ^{2}\right)  . \label{q1}%
\end{equation}

\end{theorem}

\begin{proof}
According to theorem \ref{ThConjugate}, the function $f^{-1}W_{1}$ is a
solution of (\ref{divf2}). From (\ref{reldivgrad}) we have that
\[
p^{-1/2}\left(  \operatorname{div}p\operatorname{grad}+q\right)
(p^{-1/2}W_{1})=f^{-1}\operatorname*{div}(f^{2}\nabla\left(  f^{-1}%
W_{1}\right)  )
\]
from which we obtain that $u=p^{-1/2}W_{1}$ is a solution of (\ref{maineq}).

In order to obtain the second assertion of the theorem, let us show that%
\[
p^{1/2}(\operatorname*{div}\frac{1}{p}\operatorname*{grad}+q_{1}%
)(p^{1/2}\varphi)=f\operatorname*{div}(f^{-2}\nabla(f\varphi))
\]
for any real valued $\varphi\in C^{2}(\Omega)$. According to (\ref{factUhlm}),%
\[
f\operatorname*{div}(f^{-2}\nabla(f\varphi))=\left(  \Delta-\frac{\Delta
f^{-1}}{f^{-1}}\right)  \varphi=\left(  \Delta-r_{2}\right)  \varphi.
\]
Straightforward calculation gives us the following equality%
\[
\frac{\Delta f^{-1}}{f^{-1}}=\frac{3}{4}\left(  \frac{\nabla p}{p}\right)
^{2}-\frac{1}{2}\frac{\Delta p}{p}+\left\langle \frac{\nabla p}{p}%
,\frac{\nabla u_{0}}{u_{0}}\right\rangle -\frac{\Delta u_{0}}{u_{0}}+2\left(
\frac{\nabla u_{0}}{u_{0}}\right)  ^{2}.
\]
From the condition that $u_{0}$ is a solution of (\ref{maineq}) we obtain the
equality%
\[
-\frac{\Delta u_{0}}{u_{0}}=\frac{q}{p}+\left\langle \frac{\nabla p}{p}%
,\frac{\nabla u_{0}}{u_{0}}\right\rangle .
\]
Thus,
\[
\frac{\Delta f^{-1}}{f^{-1}}=\frac{3}{4}\left(  \frac{\nabla p}{p}\right)
^{2}-\frac{1}{2}\frac{\Delta p}{p}+2\left\langle \frac{\nabla p}{p}%
,\frac{\nabla u_{0}}{u_{0}}\right\rangle +\frac{q}{p}+2\left(  \frac{\nabla
u_{0}}{u_{0}}\right)  ^{2}.
\]
Notice that
\[
\frac{\Delta p^{-1/2}}{p^{-1/2}}=\frac{3}{4}\left(  \frac{\nabla p}{p}\right)
^{2}-\frac{1}{2}\frac{\Delta p}{p}.
\]
Then
\[
\frac{\Delta f^{-1}}{f^{-1}}=\frac{\Delta p^{-1/2}}{p^{-1/2}}+2\left\langle
\frac{\nabla p}{p},\frac{\nabla u_{0}}{u_{0}}\right\rangle +\frac{q}%
{p}+2\left(  \frac{\nabla u_{0}}{u_{0}}\right)  ^{2}.
\]
Now taking $q_{1}$ in the form (\ref{q1}) we obtain the result from
(\ref{firstfact}).
\end{proof}

\subsection{Conjugate metaharmonic functions}

Theorems \ref{ThConjugate} and \ref{ThConjugate2} show us that as much as real
and imaginary parts of a complex analytic function are harmonic functions, the
real and imaginary parts of a solution of the main Vekua equation
(\ref{Vekuamain}) are solutions of associated stationary Schr\"{o}dinger
equations being also related to conductivity equations as well as to more
general elliptic equations (\ref{maineq}) and (\ref{assocmaineq}). The
following natural question arises then. We know that given an arbitrary real
valued harmonic function in a simply connected domain, a conjugate harmonic
function can be constructed explicitly such that the obtained couple of
harmonic functions represent the real and imaginary parts of a complex
analytic function. What is the corresponding more general fact for solutions
of associated stationary Schr\"{o}dinger equations (which we slightly
generalizing the definition of I. N. Vekua call metaharmonic functions) and of
other aforementioned elliptic equations. The precise result for the
Schr\"{o}dinger equations is given in the following theorem.

\begin{theorem}
\label{PrTransform}\cite{Krpseudoan} Let $W_{1}$ be a real valued solution of
(\ref{Schr1}) in a simply connected domain $\Omega$. Then the real valued
function $W_{2},$ solution of (\ref{Schr2}) such that $W=W_{1}+iW_{2}$ is a
solution of (\ref{Vekuamain}), is constructed according to the formula%
\begin{equation}
W_{2}=f^{-1}\overline{A}(if^{2}\partial_{\overline{z}}(f^{-1}W_{1})).
\label{transfDarboux}%
\end{equation}

Given a solution $W_{2}$ of (\ref{Schr2}), the corresponding solution $W_{1}$
of (\ref{Schr1}) such that $W=W_{1}+iW_{2}$ is a solution of (\ref{Vekuamain}%
), is constructed as follows%
\begin{equation}
W_{1}=-f\overline{A}(if^{-2}\partial_{\overline{z}}(fW_{2})).
\label{transfDarbouxinv}%
\end{equation}

\end{theorem}

\begin{proof}
Consider equation (\ref{Vekuamain}). Let $W=\phi f+i\psi/f$ be its solution.
Then the equation
\begin{equation}
\psi_{\overline{z}}-if^{2}\phi_{\overline{z}}=0 \label{seckind}%
\end{equation}
is valid. Note that if $W_{1}=\operatorname{Re}W$ then $\phi=W_{1}/f$. Given
$\phi$, $\psi$ is easily found from (\ref{seckind}):%
\[
\psi=\overline{A}(if^{2}\phi_{\overline{z}}).
\]
It can be verified that the expression $\overline{A}(if^{2}\phi_{\overline{z}%
})$ makes sense, that is $\partial_{x}(f^{2}\phi_{x})+\partial_{y}(f^{2}%
\phi_{y})=0.$

By theorem \ref{ThConjugate} the function $W_{2}=f^{-1}\psi$ is a solution of
(\ref{Schr2}). Thus we obtain (\ref{transfDarboux}). Let us notice that as the
operator $\overline{A}$ reconstructs the real function up to an arbitrary real
constant, the function $W_{2}$ in the formula (\ref{transfDarboux}) is
uniquely determined up to an additive term $cf^{-1}$ where $c$ is an arbitrary
real constant.

Equation (\ref{transfDarbouxinv}) is proved in a similar way.
\end{proof}

\begin{remark}
When in (\ref{Schr1}) $r_{1}\equiv0$ and $f\equiv1$, equalities
(\ref{transfDarboux}) and (\ref{transfDarbouxinv}) turn into the well known
formulas in complex analysis for constructing conjugate harmonic functions.
\end{remark}

\begin{corollary}
\cite{KrJPhys06} \label{ThDarbouxstatic} Let $U$ be a solution of
(\ref{divf2}). Then a solution $V$ of (\ref{divf-2}) such that
\[
W=fU+if^{-1}V
\]
is a solution of (\ref{Vekuamain}), is constructed according to the formula%
\begin{equation}
V=\overline{A}(if^{2}U_{\overline{z}}). \label{transfDarbouxConductivity}%
\end{equation}
Conversely, given a solution $V$ of (\ref{divf-2}), the corresponding solution
$U$ of (\ref{divf2}) can be constructed as follows:%
\[
U=-\overline{A}(if^{-2}V_{\overline{z}}).
\]

\end{corollary}

\begin{proof}
Consists in substitution of $W_{1}=fU$ and of $W_{2}=f^{-1}V$ into
(\ref{transfDarboux}) and (\ref{transfDarbouxinv}).\bigskip
\end{proof}

\begin{corollary}
\cite{KrJPhys06} \label{CorConjugate}Let $f=p^{1/2}u_{0}$, where $u_{0}$ is a
positive solution of (\ref{maineq}) in a simply connected domain $\Omega$ and
$u$ be a solution of (\ref{maineq}). Then a solution $v$ of (\ref{assocmaineq}%
) with $q_{1}$ defined by (\ref{q1}) such that $W=p^{1/2}u+ip^{-1/2}v$ is a
solution of (\ref{Vekuamain}), is constructed according to the formula%
\[
v=u_{0}^{-1}\overline{A}(ipu_{0}^{2}\partial_{\overline{z}}(u_{0}^{-1}u)).
\]

Let $v$ be a solution of (\ref{assocmaineq}), then the corresponding solution
$u$ of (\ref{maineq}) such that $W=p^{1/2}u+ip^{-1/2}v$ is a solution of
(\ref{Vekuamain}), is constructed according to the formula%
\[
u=-u_{0}\overline{A}(ip^{-1}u_{0}^{-2}\partial_{\overline{z}}(u_{0}v)).
\]

\end{corollary}

\begin{proof}
Consists in substitution of $f=p^{1/2}u_{0}$, $W_{1}=p^{1/2}u$ and
$W_{2}=p^{-1/2}v$ into (\ref{transfDarboux}) and (\ref{transfDarbouxinv}%
).\bigskip
\end{proof}

\subsection{The main Vekua equation\label{SectTheMainVekua}}

The results of this section show us that the theory of the elliptic equation
\[
(\operatorname{div}p\operatorname{grad}+q)u=0
\]
is closely related to the equation (\ref{Vekuamain}):
\begin{equation}
W_{\overline{z}}=\frac{f_{\overline{z}}}{f}\overline{W}.\label{Vekuamain1}%
\end{equation}
As was pointed out in remark \ref{RelAstalaPaivarinta} the pair of functions:
$F=f\quad$and\quad$G=\frac{i}{f}$ is a generating pair for this equation. Then
the corresponding characteristic coefficients $A_{(F,G)}$ and $B_{(F,G)}$ have
the form%
\[
A_{(F,G)}=0,\quad\text{\quad}B_{(F,G)}=\frac{f_{z}}{f},
\]
and the $(F,G)$-derivative according to (\ref{FGder}) is defined as follows%
\[
\overset{\cdot}{W}=W_{z}-\frac{f_{z}}{f}\overline{W}=\left(  \partial
_{z}-\frac{f_{z}}{f}C\right)  W.
\]

Comparing $B_{(F,G)}$ with the coefficient in (\ref{Vekua}) and due to Theorem
\ref{ThBersDer} we obtain the following statement.

\begin{proposition}
\label{PrDer} Let $W$ be a solution of (\ref{Vekuamain1}). Then its
$(F,G)$-derivative, the function $w=\overset{\cdot}{W}$ is a solution of
(\ref{Vekua}).
\end{proposition}

This result can be verified also by a direct substitution.

According to (\ref{FGAnt}) and taking into account that
\[
F^{\ast}=-if\quad\text{and\quad}G^{\ast}=1/f,
\]
the $(F,G)$-antiderivative has the form%
\begin{align}
\int_{z_{0}}^{z}w(\zeta)d_{(F,G)}\zeta &  =f(z)\operatorname{Re}\int_{z_{0}%
}^{z}\frac{w(\zeta)}{f(\zeta)}d\zeta-\frac{i}{f(z)}\operatorname{Re}%
\int_{z_{0}}^{z}if(\zeta)w(\zeta)d\zeta\nonumber\\
&  =f(z)\operatorname{Re}\int_{z_{0}}^{z}\frac{w(\zeta)}{f(\zeta)}d\zeta
+\frac{i}{f(z)}\operatorname{Im}\int_{z_{0}}^{z}f(\zeta)w(\zeta)d\zeta,
\label{antider}%
\end{align}
and we obtain the following statement.

\begin{proposition}
\label{PrAntider} Let $w$ be a solution of (\ref{Vekua}). Then the function
\[
W(z)=\int_{z_{0}}^{z}w(\zeta)d_{(F,G)}\zeta
\]
is a solution of (\ref{Vekuamain1}).
\end{proposition}

\subsection{$p$-analytic functions}

\begin{definition}
A function $\Phi=u+iv$ of a complex variable $z=x+iy$ is said to be
$p$-analytic in some domain $\Omega$ iff
\begin{equation}
u_{x}=\frac{1}{p}v_{y},\qquad u_{y}=-\frac{1}{p}v_{x}\qquad\text{in }%
\Omega\label{Polozh}%
\end{equation}
where $p$ is a given positive function of $x$ and $y$ which is supposed to be
continuously differentiable.
\end{definition}

The theory of $p$-analytic functions was presented in \cite{Polozhy}.
$p$-analytic functions in a certain sense represent a subclass of generalized
analytic (or pseudoanalytic) functions studied by L. Bers \cite{Berskniga,
BersStat} and I. N. Vekua \cite{Vekua} and it should be noticed that this
subclass preserves some important properties of usual analytic functions which
are not preserved by a too ample class of generalized analytic functions
(corresponding details can be found in \cite{Polozhy}). $p$-analytic functions
found applications in elasticity theory (see, e.g., \cite{AS, Goman}) and in
axisymmetric problems of hydrodynamics (see, e.g., \cite{Garab},
\cite{ZabarankinUlitko2006}).

As was pointed out in remark \ref{RelAstalaPaivarinta} the function $W=\phi
f+i\psi/f$ is a solution of the main Vekua equation (\ref{Vekuamain1}) if and
only if $\phi$ and $\psi$ satisfy the equation (\ref{Vekuamain2kind}) which is
equivalent to the system%
\[
\phi_{x}=\frac{1}{f^{2}}\psi_{y},\qquad\phi_{y}=-\frac{1}{f^{2}}\psi_{x}.
\]
In other words $W$ is a solution of the main Vekua equation iff its
corresponding pseudoanalytic function of second kind is $f^{2}$-analytic.
Thus, we obtain the following connection between the stationary
Schr\"{o}dinger equation and the system defining $p$-analytic functions.

\begin{theorem}
Let $f$ be a positive solution of the equation%
\begin{equation}
-\Delta g+\nu g=0 \label{Schrod1A}%
\end{equation}
where $\nu$ is a real valued function and let $W_{1}$ be another real valued
solution of this equation. Then the function $\Phi=W_{1}/f+ifW_{2}$, where
$W_{2}$ is defined by (\ref{transfDarboux}) is an $f^{2}$-analytic function,
and vice versa, let $\Phi$ be an $f^{2}$-analytic function then the function
$W_{1}=f\operatorname*{Re}\Phi$ is a solution of (\ref{Schrod1A}).
\end{theorem}

The following relation between solutions of the conductivity equation and
$p$-analytic functions is valid also.

\begin{theorem}
Let $f$ be a positive continuously differentiable function in a domain
$\Omega$ and let $U$ be a real valued solution of the equation
\begin{equation}
\operatorname*{div}(f^{2}\nabla U)=0\qquad\text{in }\Omega.
\label{conductivity}%
\end{equation}
Then the function $\Phi=U+iV$ is $f^{2}$-analytic in $\Omega$, where $V$ is
defined by (\ref{transfDarbouxConductivity}), and vice versa, let $\Phi$ be
$f^{2}$-analytic in $\Omega$ then $U=\operatorname*{Re}\Phi$ is a solution of
(\ref{conductivity}).
\end{theorem}

Thus, solutions of the stationary Schr\"{o}dinger equation and of the
conductivity equation can be converted into $p$-analytic functions and vice
versa. In some cases this relation leads to a simplification of a part of
$p$-analytic function theory.

\begin{example}
\cite{KrPanalyt} A considerable part of bibliography dedicated to $p$-analytic
functions consists of studying the case $p=x^{k}$, where $k\in\mathbf{R}$
(see, e.g., \cite{Chemeris}, \cite{KapshKlen}, \cite{Polozhy}). Let us see
what is the form of the corresponding Schr\"{o}dinger equation. For this we
should calculate the potential $\nu$ in (\ref{Schrod1A}) when $f=x^{k/2}$. It
is easy to see that
\begin{equation}
\nu=\frac{k^{2}-2k}{4x^{2}}. \label{nu}%
\end{equation}
The Schr\"{o}dinger equation with this potential is well studied. Separation
of variables leads us to the equation%
\begin{equation}
X^{\prime\prime}(x)+\left(  \beta^{2}-\frac{4\alpha^{2}-1}{4x^{2}}\right)
X(x)=0, \label{Bessel}%
\end{equation}
where $\beta^{2}$ is the separation constant and $\alpha=(k-1)/2$. The
function
\[
X(x)=\sqrt{x}Z_{\alpha}(\beta x)
\]
is a solution of (\ref{Bessel}) (see \cite[8.491]{GradRizh}) where $Z_{\alpha
}$ denotes any cylindric function of order $\alpha$ (Bessel functions of first
or second kind). Thus the study of $x^{k}$-analytic functions reduces to the
Schr\"{o}dinger equation (\ref{Schrod}) with $\nu$ defined by (\ref{nu}) which
in its turn after having separated variables reduces to a kind of Bessel
equation (\ref{Bessel}).
\end{example}

\begin{example}
In the work \cite{KapshYazk} boundary value problems for $p$-analytic
functions with $p=x/(x^{2}+y^{2})$ were studied. Considering
\[
f=\sqrt{p}=\sqrt{\frac{x}{x^{2}+y^{2}}}%
\]
we see that this function is a solution of the Schr\"{o}dinger equation
(\ref{Schrod}) with $\nu$ having the form%
\[
\nu=-\frac{1}{4x^{2}},
\]
that is we obtain again the potential of the form (\ref{nu}) where $k=1$ and
as was shown in the previous Remark the study of corresponding $p$-analytic
functions in a sense reduces to the Bessel equation (\ref{Bessel}).
\end{example}

\section{Complete systems of solutions for second-order
equations\label{SectComplSystems}}

In what follows let us suppose that the real valued function $f$ is defined in
a somewhat bigger domain $\Omega_{\varepsilon}$ with a sufficiently smooth
boundary. Then we change the function $f$ for $z\in\Omega_{\varepsilon}^{{}%
}\backslash\Omega$ and continue it over the whole plane in such a way that
$f\equiv1$ for large $\left\vert z\right\vert $ (see \cite{BersFormalPowers}).
In this way the generating pair $(F,G)=(f,i/f)$ becomes complete and normalized.

Then the following statements are direct corollaries of the relations
established in section \ref{SectSchr} between pseudoanalytic functions
(solutions of (\ref{Vekuamain})) and solutions of second-order elliptic
equations, and of the convergence theorems from subsection
\ref{SubsectConvergence}.

\begin{definition}
Let $u(z)$ be a given solution of the equation (\ref{maineq}) defined for
small values of $\left\vert z-z_{0}\right\vert $, and let $W(z)$ be a solution
of (\ref{Vekuamain}) constructed according corollary \ref{CorConjugate}, such
that $\operatorname*{Re}W=p^{1/2}u$. The series
\[
p^{-1/2}(z)\sum_{n=0}^{\infty}\operatorname*{Re}Z^{(n)}(a_{n},z_{0};z)
\]
with the coefficients given by (\ref{Taylorcoef}) is called the Taylor series
of $u$ at $z_{0}$, formed with formal powers.
\end{definition}

\begin{theorem}
Let $u(z)$ be a solution of (\ref{maineq}) defined for $\left\vert
z-z_{0}\right\vert <R$. Then it admits a unique expansion of the form
\[
u(z)=p^{-1/2}(z)\sum_{n=0}^{\infty}\operatorname*{Re}Z^{(n)}(a_{n},z_{0};z)
\]
which converges normally for $\left\vert z-z_{0}\right\vert <R$.
\end{theorem}

\begin{proof}
This is a direct consequence of theorem \ref{ThTaylorRepr} and remark
\ref{RemNecConditions}. Both necessary conditions in remark
\ref{RemNecConditions} are fulfilled for the generating pair (\ref{genpair}).
\end{proof}

\begin{theorem}
\label{ThRungeSchr}An arbitrary solution of (\ref{maineq}) defined in a simply
connected domain where there exists a positive particular solution $u_{0}$ can
be expanded into a normally convergent series of formal polynomials multiplied
by $p^{-1/2}$.
\end{theorem}

\begin{proof}
This is a direct corollary of theorem \ref{ThRunge}.
\end{proof}

More precisely the last theorem has the following meaning. Due to Property 2
of formal powers we have that $Z^{(n)}(a,z_{0};z)$ for any Taylor coefficient
$a$ can be expressed through $Z^{(n)}(1,z_{0};z)$ and $Z^{(n)}(i,z_{0};z)$.
Then due to theorem \ref{ThRunge} any solution $W$ of (\ref{Vekuamain}) can be
expanded into a normally convergent series of linear combinations of
$Z^{(n)}(1,z_{0};z)$ and $Z^{(n)}(i,z_{0};z)$. Consequently, any solution of
(\ref{maineq}) can be expanded into a normally convergent series of linear
combinations of real parts of $Z^{(n)}(1,z_{0};z)$ and $Z^{(n)}(i,z_{0};z)$
multiplied by $p^{-1/2}$.

Obviously, for solutions of (\ref{maineq}) the results on the interpolation
and on the degree of approximation like, e.g., theorem \ref{ThMenke} are also valid.

Let us stress that theorem \ref{ThRungeSchr} gives us the following result.
The functions
\begin{equation}
\left\{  p^{-1/2}(z)\operatorname*{Re}Z^{(n)}(1,z_{0};z),\quad p^{-1/2}%
(z)\operatorname*{Re}Z^{(n)}(i,z_{0};z)\right\}  _{n=0}^{\infty}%
\label{complsystem}%
\end{equation}
represent a complete system of solutions of (\ref{maineq}) in the sense that
any solution of (\ref{maineq}) can be represented by a normally convergent
series formed by functions (\ref{complsystem}) in any simply connected domain
$\Omega$ where a positive solution of (\ref{maineq}) exists. Moreover, as we
show in section 6, in many practically interesting situations these functions
can be constructed explicitly.

\section{A remark on orthogonal coordinate systems in a plane}

Orthogonal coordinate systems in a plane are obtained (see \cite{Madelung})
from Cartesian coordinates $x$, $y$ by means of the relation%
\[
u+iv=\Phi(x+iy)
\]
where $\Phi$ is an arbitrary analytic function. Quite often a transition to
more general coordinates is useful
\[
\xi=\xi(u),\quad\eta=\eta(v).
\]
$\xi$ and $\eta$ preserve the property of orthogonality. Some examples taken
from \cite{Madelung} illustrate the point.

\begin{example}
\textbf{Polar coordinates }%
\[
u+iv=\ln(x+iy),
\]%
\begin{equation}
u=\ln\sqrt{x^{2}+y^{2}},\quad v=\arctan\frac{y}{x}. \label{trueorthogcoord}%
\end{equation}
Usually the following new coordinates are introduced%
\[
r=e^{u}=\sqrt{x^{2}+y^{2}},\quad\varphi=v=\arctan\frac{y}{x}.
\]

\end{example}

\begin{example}
\textbf{Parabolic coordinates }%
\[
\frac{u+iv}{\sqrt{2}}=\sqrt{x+iy},
\]%
\[
u=\sqrt{r+x},\quad v=\sqrt{r-x}.
\]
More frequently the parabolic coordinates are introduced as follows%
\[
\xi=u^{2},\quad\eta=v^{2}.
\]

\end{example}

\begin{example}
\textbf{Elliptic coordinates }%
\[
u+iv=\arcsin\frac{x+iy}{\alpha},
\]%
\[
\sin u=\frac{s_{1}-s_{2}}{2\alpha},\quad\cosh v=\frac{s_{1}+s_{2}}{2\alpha}%
\]
where $s_{1}=\sqrt{(x+\alpha)^{2}+y^{2}}$, $s_{2}=\sqrt{(x-\alpha)^{2}+y^{2}}%
$. The substitution
\[
\xi=\sin u,\quad\eta=\cosh v
\]
is frequently used.
\end{example}

\begin{example}
\textbf{Bipolar coordinates}%
\[
u+iv=\ln\frac{\alpha+x+iy}{\alpha-x-iy},
\]%
\[
\tanh u=\frac{2\alpha x}{\alpha^{2}+x^{2}+y^{2}},\quad\tan v=\frac{2\alpha
y}{\alpha^{2}-x^{2}-y^{2}}.
\]
The following substitution is frequently used%
\[
\xi=e^{-u},\quad\eta=\pi-v.
\]

\end{example}

\section{Explicit construction of a generating
sequence\label{SectExplicitConstrGenSeq}}

We suppose that the function $f$ in the main Vekua equation (\ref{Vekuamain})
has the following form%
\begin{equation}
f=U(u)V(v)\label{fUV}%
\end{equation}
where $u$ and $v$ represent an orthogonal coordinate system and according to
the explained in the previous section we assume that $\Phi=u+iv$ is an
analytic function of the variable $z=x+iy$. $U$ and $V$ are arbitrary
differentiable nonvanishing real valued functions.

The present section is structured in the following way. First we explain how
one can naturally arrive at the form of generating sequences in the main
result of the section, theorem \ref{ThGenSeq} and then we give a rigorous
proof of this result.

The first step in the construction of a generating sequence for the main Vekua
equation (\ref{Vekuamain}) is the construction of a generating pair for the
equation (\ref{Vekua}) which as was shown in subsection \ref{SectTheMainVekua}
is a successor of the main Vekua equation. For this one of the possibilities
consists in constructing another pair of solutions of (\ref{Vekuamain}). Then
their $(F,G)$-derivatives will give us solutions of (\ref{Vekua}). Let us
consider equation (\ref{Vekuamain}) in its equivalent form
(\ref{Vekuamain2kind}):
\begin{equation}
\varphi_{\overline{z}}+\frac{i}{f^{2}}\psi_{\overline{z}}=0 \label{seckindA}%
\end{equation}
and look for a solution in the form%
\[
\varphi=\varphi(u)\text{,}\qquad\psi=\psi(v).
\]
Then we have the following equation
\[
\varphi^{\prime}(u)u_{\overline{z}}+\frac{i}{f^{2}}\psi^{\prime}%
(v)v_{\overline{z}}=0.
\]
Taking into account that $\Phi=u+iv$ is an analytic function, that is
$u_{\overline{z}}+iv_{\overline{z}}=0$ we observe that (\ref{seckindA}) is
fulfilled if $\varphi^{\prime}(u)=\psi^{\prime}(v)/f^{2}.$ Now using
(\ref{fUV}) we obtain that $U^{2}(u)V^{2}(v)=\psi^{\prime}(v)/\varphi^{\prime
}(u)$ and hence%
\[
\varphi(u)=\int\frac{du}{U^{2}(u)}\text{\quad and\quad}\psi(v)=\int
V^{2}(v)dv.
\]
The corresponding solution $W_{1}$ of (\ref{Vekuamain}) has the form%
\[
W_{1}=\int\frac{du}{U^{2}(u)}U(u)V(v)+\int V^{2}(v)dv\frac{i}{U(u)V(v)}.
\]
Its $(F,G)$-derivative is obtained according to (\ref{FGder}) as follows%
\[
\overset{\cdot}{W}_{1}=\frac{V}{U}u_{z}+i\frac{V}{U}v_{z}=\frac{V}{U}\Phi
_{z}.
\]

By analogy, we can look for a solution of (\ref{seckindA}) in the form
\[
\varphi=\varphi(v)\text{,}\qquad\psi=\psi(u).
\]
Then we have the following equation
\[
\varphi^{\prime}(v)v_{\overline{z}}+\frac{i}{f^{2}}\psi^{\prime}%
(u)u_{\overline{z}}=0.
\]
This equation is fulfilled if $\varphi^{\prime}(v)=-1/V^{2}(v)$ and
$\psi^{\prime}(u)=U^{2}(u).$ Consequently,%
\[
\varphi(v)=-\int\frac{dv}{V^{2}(v)}\text{\quad and\quad}\psi(u)=\int
U^{2}(u)du.
\]
The corresponding solution $W_{2}$ of (\ref{Vekuamain}) has the form%
\[
W_{2}=-\int\frac{dv}{V^{2}(v)}U(u)V(v)+\int U^{2}(u)du\frac{i}{U(u)V(v)}%
\]
with its corresponding $(F,G)$-derivative
\[
\overset{\cdot}{W}_{2}=-\frac{U}{V}v_{z}+i\frac{U}{V}u_{z}=i\frac{U}{V}%
\Phi_{z}.
\]
Denote $F_{1}=\overset{\cdot}{W}_{1}$ and $G_{1}=\overset{\cdot}{W}_{2}$. It
is easy to see that the pair $F_{1},$ $G_{1}$ fulfils (\ref{GenPairCond}) and
by construction satisfies equation (\ref{Vekua}). Thus, $(F_{1},G_{1})$ is a
successor of $(F,G)$ defined by (\ref{genpair}).

The following step is to construct the generating pair $(F_{2},G_{2})$. For
this we should find another pair of solutions of (\ref{Vekua}) and apply the
$(F_{1},G_{1})$-derivative to them. Consider equation (\ref{Vekua}) written in
the form%
\[
\varphi_{\overline{z}}F_{1}+\psi_{\overline{z}}G_{1}=0
\]
which in our case can be represented as follows%
\begin{equation}
\varphi_{\overline{z}}\frac{V}{U}+\psi_{\overline{z}}i\frac{U}{V}=0.
\label{Vekua1}%
\end{equation}
Again, let us look for a solution in the form%
\[
\varphi=\varphi(u)\text{,}\qquad\psi=\psi(v).
\]
Then equation (\ref{Vekua1}) is satisfied if $\varphi^{\prime}(u)=\left(
\frac{U(u)}{V(v)}\right)  ^{2}\psi^{\prime}(v)$ from where we obtain
\[
\varphi(u)=\int U^{2}(u)du\text{\quad and\quad}\psi(v)=\int V^{2}(v)dv.
\]
Then the corresponding solution of (\ref{Vekua}) has the form
\[
w_{1}=\int U^{2}(u)du\Phi_{z}\frac{V(v)}{U(u)}+\int V^{2}(v)dv\Phi_{z}%
\frac{iU(u)}{V(v)}.
\]
Its $(F_{1},G_{1})$-derivative is obtained as follows
\[
\overset{\cdot}{w}_{1}=\Phi_{z}UVu_{z}+i\Phi_{z}UVv_{z}=UV\left(  \Phi
_{z}\right)  ^{2}.
\]
Analogously, looking for a solution of (\ref{Vekua1}) in the form
\[
\varphi=\varphi(v)\text{,}\qquad\psi=\psi(u)
\]
we obtain that
\[
\varphi(v)=-\int\frac{dv}{V^{2}(v)}\text{\quad and\quad}\psi(u)=\int\frac
{du}{U^{2}(u)}.
\]
The corresponding solution of (\ref{Vekua}) has the form
\[
w_{2}=-\int\frac{dv}{V^{2}(v)}\Phi_{z}\frac{V(v)}{U(u)}+\int\frac{du}%
{U^{2}(u)}\Phi_{z}\frac{iU(u)}{V(v)}.
\]
Its $(F_{1},G_{1})$-derivative is obtained as follows%
\[
\overset{\cdot}{w}_{2}=-\frac{\Phi_{z}v_{z}}{UV}+\frac{\Phi_{z}iu_{z}}%
{UV}=\frac{i}{UV}\left(  \Phi_{z}\right)  ^{2}.
\]
We have that the pair $F_{2}=\overset{\cdot}{w}_{1}$ and $G_{2}=\overset
{\cdot}{w}_{2}$is a successor of $(F_{1},G_{1})$. Observe that
\[
(F_{1},G_{1})=\left(  \frac{\Phi_{z}}{U^{2}}F,\quad U^{2}\Phi_{z}G\right)
\]
and%
\[
(F_{2},G_{2})=\left(  \left(  \Phi_{z}\right)  ^{2}F,\quad\left(  \Phi
_{z}\right)  ^{2}G\right)  .
\]
From these formulae it is easy to guess a general form of the corresponding
generating sequence which we give in the following statement.

\begin{theorem}
\label{ThGenSeq} Let $F=U(u)V(v)$ and $G=\frac{i}{U(u)V(v)}$ where $U$ and $V$
are arbitrary differentiable nonvanishing real valued functions, $\Phi=u+iv$
is an analytic function of the variable $z=x+iy$ in $\Omega$ such that
$\Phi_{z}$ is bounded and has no zeros in $\Omega$. Then the generating pair
$(F,G)$ is embedded in the generating sequence $(F_{m},G_{m})$, $m=0,\pm
1,\pm2,\ldots$in $\Omega$ defined as follows
\[
F_{m}=\left(  \Phi_{z}\right)  ^{m}F\quad\text{and\quad}G_{m}=\left(  \Phi
_{z}\right)  ^{m}G\quad\text{for even }m
\]
and%
\[
F_{m}=\frac{\left(  \Phi_{z}\right)  ^{m}}{U^{2}}F\quad\text{and\quad}%
G_{m}=\left(  \Phi_{z}\right)  ^{m}U^{2}G\quad\text{for odd }m.
\]

\end{theorem}

\begin{proof}
First of all let us show that $(F_{m},G_{m})$ is a generating pair for
$m=\pm1,\pm2,\ldots$ . Indeed we have
\[
\operatorname{Im}(\overline{F}_{m}G_{m})=\operatorname{Im}(\left\vert \Phi
_{z}\right\vert ^{2m}\overline{F}G)>0.
\]

Consider $a_{(F_{m},G_{m})}$. For both $m$ being even or odd we obtain%
\[
a_{(F_{m},G_{m})}=\left\vert \Phi_{z}\right\vert ^{2m}a_{(F,G)}\equiv0.
\]
We should verify the equality
\begin{equation}
b_{(F_{m},G_{m})}=-B_{(F_{m-1},G_{m-1})} \label{vsp}%
\end{equation}
Consider first the case of an odd $m$. A direct calculation gives us%
\[
b_{(F_{m},G_{m})}=\left(  \frac{\Phi_{z}^{{}}}{\overline{\Phi_{z}}}\right)
^{m}\left(  b_{(F,G)}-2u_{\overline{z}}\frac{U^{\prime}}{U}\right)
\]
and
\[
B_{(F_{m-1},G_{m-1})}=\left(  \frac{\Phi_{z}^{{}}}{\overline{\Phi_{z}}%
}\right)  ^{m-1}B_{(F,G)}.
\]
Thus equality (\ref{vsp}) is true iff $\frac{\Phi_{z}^{{}}}{\overline{\Phi
_{z}}}\left(  b_{(F,G)}-2u_{\overline{z}}\frac{U^{\prime}}{U}\right)  $ is
equal to $-B_{(F,G)}$. It is easy to see that
\[
B_{(F,G)}=u_{z}\left(  \frac{U^{\prime}}{U}-i\frac{V^{\prime}}{V}\right)  .
\]
Consider%
\[
\frac{\Phi_{z}^{{}}}{\overline{\Phi_{z}}}\left(  b_{(F,G)}-2u_{\overline{z}%
}\frac{U^{\prime}}{U}\right)  =\frac{u_{z}+iv_{z}}{u_{\overline{z}%
}-iv_{\overline{z}}}\left(  \frac{U^{\prime}}{U}u_{\overline{z}}%
+\frac{V^{\prime}}{V}v_{\overline{z}}-2u_{\overline{z}}\frac{U^{\prime}}%
{U}\right)
\]%
\[
=\frac{u_{z}}{u_{\overline{z}}}\left(  -\frac{U^{\prime}}{U}u_{\overline{z}%
}+\frac{V^{\prime}}{V}iu_{\overline{z}}\right)  =u_{z}\left(  -\frac
{U^{\prime}}{U}+i\frac{V^{\prime}}{V}\right)  .
\]
Thus, equality (\ref{vsp}) is proved in the case of $m$ being odd.

Now let $m$ be even. Then
\[
b_{(F_{m},G_{m})}=\left(  \frac{\Phi_{z}^{{}}}{\overline{\Phi_{z}}}\right)
^{m}b_{(F,G)}%
\]
and%
\[
B_{(F_{m-1},G_{m-1})}=\left(  \frac{\Phi_{z}^{{}}}{\overline{\Phi_{z}}%
}\right)  ^{m-1}\left(  B_{(F,G)}-2\frac{U^{\prime}}{U}u_{z}\right)  .
\]
Equality (\ref{vsp}) is valid iff the expression $\frac{\Phi_{z}^{{}}%
}{\overline{\Phi_{z}}}b_{(F,G)}$ is equal to $\left(  -B_{(F,G)}%
+2\frac{U^{\prime}}{U}u_{z}\right)  $. We have
\[
\frac{\Phi_{z}^{{}}}{\overline{\Phi_{z}}}b_{(F,G)}=\frac{u_{z}+iv_{z}%
}{u_{\overline{z}}-iv_{\overline{z}}}\left(  \frac{U^{\prime}}{U}%
u_{\overline{z}}+\frac{V^{\prime}}{V}v_{\overline{z}}\right)  =\frac{u_{z}%
}{u_{\overline{z}}}\left(  \frac{U^{\prime}}{U}u_{\overline{z}}+\frac
{V^{\prime}}{V}iu_{\overline{z}}\right)  =u_{z}\left(  \frac{U^{\prime}}%
{U}+i\frac{V^{\prime}}{V}\right)
\]
and from the other side
\[
-B_{(F,G)}+2\frac{U^{\prime}}{U}u_{z}=-u_{z}\left(  \frac{U^{\prime}}%
{U}-i\frac{V^{\prime}}{V}\right)  +2\frac{U^{\prime}}{U}u_{z}=u_{z}\left(
\frac{U^{\prime}}{U}+i\frac{V^{\prime}}{V}\right)  .
\]
Thus, equality (\ref{vsp}) is proved in all cases and the sequence
$(F_{m},G_{m})$, $m=0,\pm1,\pm2,\ldots$ satisfies the conditions of Definition
\ref{DefSeq}. Therefore it is a generating sequence.
\end{proof}

\begin{remark}
This result obviously generalizes the explicit construction of a generating
sequence in the case when $u=x$ and $v=y$ presented in a number of works by L.
Bers (see, e.g., \cite{Berskniga} and \cite{Courant}).
\end{remark}

The last theorem opens the way for explicit construction of formal powers
corresponding to the main Vekua equation (\ref{Vekuamain}) in the case when
$f$ has the form (\ref{fUV}) and hence for explicit construction of complete
systems of solutions for corresponding second-order elliptic equations. We
give more details on this as well as some examples in the next section.

\section{Explicit construction of complete systems of solutions of
second-order elliptic equations}

In section \ref{SectComplSystems} we explained how complete systems of
solutions of second-order elliptic equations can be constructed from systems
of formal powers for a corresponding main Vekua equation (\ref{Vekuamain}). In
the preceding section we established a result which allows us to make this
construction possible in the case when $f$ in (\ref{Vekuamain}) has the form
(\ref{fUV}). Then formal powers are constructed simply by definition
\ref{DefFormalPower}. The meaning of this result for the stationary
Schr\"{o}dinger equation and for the conductivity equation we discuss
separately in the next two subsections.

\subsection{Explicit construction of complete systems of solutions for a
stationary Schr\"{o}dinger equation}

Consider the equation
\begin{equation}
-\Delta g+\nu g=0\qquad\text{in }\Omega. \label{Schrod1B}%
\end{equation}
where $\nu$ is a real valued function. In order to start the procedure of
construction of a complete system of solutions of this equation we need a
positive particular solution $f$ having the form (\ref{fUV}).

\begin{example}
Let $\nu$ in (\ref{Schrod1B}) depend on one Cartesian variable: $\nu=\nu(x).$
Suppose we are given a particular solution $f=f(x)$ of the ordinary
differential equation
\begin{equation}
-\frac{d^{2}f}{dx^{2}}+\nu f=0. \label{Schrord}%
\end{equation}
It is sufficient for the application of our result from the preceding section
for constructing the corresponding generating sequence and hence the system of
formal powers for the main Vekua equation which in this particular case has
the form%
\[
W_{\overline{z}}=\frac{f_{x}}{2f}\overline{W}.
\]

\end{example}

\begin{example}
A number of works (see, e.g., \cite{ChadanKobayashi06A} and
\cite{ChadanKobayashi06B}) are dedicated to construction (in our terms) of a
particular solution for the Schr\"{o}dinger equation with a radially symmetric
potential. This solution $f$ is precisely the only necessary ingredient in
order to obtain a complete system of formal powers and hence of solutions of
the Schr\"{o}dinger equation. Here an important restriction is that the
analytic function $\Phi_{z}$ corresponding to the polar coordinate system has
the form $1/z$ and hence our procedure works in any domain not including the
origin where $f$ is positive.
\end{example}

\begin{example}
Consider the Yukawa equation
\begin{equation}
(-\Delta+c^{2})u=0 \label{Helmholtz}%
\end{equation}
with $c$ being a real constant. Take the following particular solution of
(\ref{Helmholtz}): $f=e^{cy}$. Let us construct the first few corresponding
formal powers with center at the origin. We have%
\[
Z^{(0)}(1,0;z)=e^{cy},\qquad Z^{(0)}(i,0;z)=ie^{-cy},
\]%
\[
Z^{(1)}(1,0;z)=xe^{cy}+\frac{i\sinh(cy)}{c},\qquad Z^{(1)}(i,0;z)=-\frac
{\sinh(cy)}{c}+ixe^{-cy},
\]%
\[
Z^{(2)}(1,0;z)=\left(  x^{2}-\frac{y}{c}\right)  e^{cy}+\frac{\sinh(cy)}%
{c^{2}}+\frac{2ix\sinh(cy)}{c},
\]%
\[
Z^{(2)}(i,0;z)=-\frac{2x\sinh(cy)}{c}+i\left(  \left(  x^{2}+\frac{y}%
{c}\right)  e^{-cy}-\frac{\sinh(cy)}{c^{2}}\right)  ,\ldots.
\]
It is a simple exercise to verify that indeed the asymptotic formulas
(\ref{asymptformulas}) hold. Now taking real parts of the formal powers we
obtain a complete system of solutions of the Yukawa equation:%
\[
u_{0}(x,y)=e^{cy},\qquad u_{1}(x,y)=xe^{cy},\qquad u_{2}(x,y)=-\frac
{\sinh(cy)}{c},
\]%
\[
u_{3}(x,y)=\left(  x^{2}-\frac{y}{c}\right)  e^{cy}+\frac{\sinh(cy)}{c^{2}%
},\qquad u_{4}(x,y)=-\frac{2x\sinh(cy)}{c},\ldots\text{.}%
\]
Formal powers of higher order can be constructed explicitly using a computer
system of symbolic calculation. For this particular example (together with
Maria Rosal\'{\i}a Tenorio) Matlab 6.5 allowed us to obtain analytic
expressions for the formal powers up to the order ten, that gave us the first
twenty one functions $u_{0},\ldots,u_{20}$. We used them for a numerical
solution of the Dirichlet problem for the Yukawa equation with very
satisfactory results. For example, in the case when $\Omega$ is a unit disk
with centre at the origin, $c=1$ and $u$ on the boundary is equal to $e^{x}$
(this test exact solution gave us the worst precision because of its obvious
\textquotedblleft disparateness\textquotedblright\ from functions $u_{0}%
,u_{1}\ldots$) the maximal error $\max_{z\in\Omega}\left\vert u(z)-\widetilde
{u}(z)\right\vert $ where $u$ is the exact solution and $\widetilde{u}%
=\sum_{n=0}^{20}a_{n}u_{n}$, the real constants $a_{n}$ being found by the
collocation method, was of order 10$^{-7}$. A very fast convergence of the
method was observed.

Although the numerical method based on the usage of explicitly or numerically
constructed pseudoanalytic formal powers still needs a much more detailed
analysis these first results show us that it is quite possible that in due
time and with a further development of symbolic calculation systems it can
rank high among other numerical approaches. The use of complete systems of
exact solutions based on pseudoanalytic formal powers has an important
advantage of universality. The system does not depend on the choice of a
domain whenever the function $f$ keeps to have no zeros.

In \cite{BraggDettman95} the following system of solutions of the Yukawa
equation (\ref{Helmholtz}) was obtained from completely different reasonings,%
\[
u_{0}=\cosh(cy),\text{\quad}u_{1}=x\cosh(cy),\quad u_{2}=x^{2}\cosh
(cy)-\frac{y}{c}\sinh(cy),
\]%
\[
u_{3}=x^{3}\cosh(cy)-\frac{3xy}{c}\sinh(cy),\ldots.
\]
The same system is easily obtained (our Matlab program gives symbolic
representations of these functions up to the subindex 21) if instead of the
particular solution $f=e^{cy}$ considered in the above example we choose
$f=\cosh(cy)$.

In the case of the Helmholtz equation%
\[
(\Delta+c^{2})u=0
\]
L. R. Bragg and J. W. Dettman constructed the following system%
\[
u_{0}=\cos(cy),\text{\quad}u_{1}=x\cos(cy),\quad u_{2}=x^{2}\cos(cy)-\frac
{y}{c}\sin(cy),
\]%
\[
u_{3}=x^{3}\cos(cy)-\frac{3xy}{c}\sin(cy),\ldots.
\]
This system is also easily obtained as explained above choosing $f=\cos(cy)$.

Thus, the construction of the systems of solutions for the Yukawa and for the
Helmholtz equation as well as the proof of their completeness in
\cite{BraggDettman95} are corollaries of a more general theory presented here.
\end{example}

In fact, having a positive particular solution $f$ of the form (\ref{fUV}) for
equation (\ref{Schrod1B}) allows us to construct explicitly a complete system
of solutions of (\ref{Schrod1B}) in any simply connected bounded domain where
the corresponding to a chosen coordinate system analytic function $\Phi_{z}$
is bounded and different from zero.

\subsection{Explicit construction of complete systems of solutions for a
conductivity equation}

Consider the equation
\[
\operatorname*{div}(f^{2}\nabla U)=0.
\]
Usually in practical applications the function $f$ is positive and equals $1$
outside a certain disk. We suppose additionally that $f$ is continuously
differentiable on the whole plane. Under these conditions $F=f$ and $G=i/f$ is
a complete normalized generating pair and all the convergence and Runge type
results are valid. Using results of section \ref{SectExplicitConstrGenSeq} one
can construct explicitly a complete system of solutions of the conductivity
equation when the function $f$ has the form (\ref{fUV}) in any domain where
the analytic function $\Phi_{z}$ from theorem \ref{ThGenSeq} is bounded and
has no zeros.

\end{document}